\documentclass[11pt]{amsart}

\usepackage{amsmath,amssymb}
\usepackage{a4wide}
\usepackage{algorithm}
\usepackage{algorithmic}

\usepackage{graphicx}

\newtheorem{lemma}{Lemma}[section]
\newtheorem{theorem}[lemma]{Theorem}

\newcommand{\dist}{\textrm{dist}}

\begin{document}


\title{Uniqueness of graph square roots of girth six}

\author[Anna Adamaszek, Micha{\l} Adamaszek]{Anna Adamaszek$^{1,3}$, Micha{\l} Adamaszek$^{2,3}$}
\thanks{Research supported by the Centre for Discrete
        Mathematics and its Applications (DIMAP), EPSRC award EP/D063191/1.}
\thanks{$^1$Department of Computer Science, University of Warwick, Coventry, CV4 7AL, UK}
\thanks{$^2$Warwick Mathematics Institute, University of Warwick, Coventry, CV4 7AL, UK}
\thanks{$^3$Centre for Discrete Mathematics and its Applications (DIMAP), University of Warwick}
\thanks{E-mails: \texttt{\{annan,aszek\}@mimuw.edu.pl}}



\begin{abstract}
We prove that if two graphs of girth at least $6$ have isomorphic squares, then the graphs themselves are isomorphic. This is the best possible extension of the results of Ross and Harary on trees and the results of Farzad et al. on graphs of girth at least $7$. We also make a remark on reconstruction of graphs from their higher powers.
\end{abstract}

\maketitle

\section{Introduction}
\label{section:introduction}

For a simple, undirected, connected graph $H$ its \emph{square} $G=H^2$ is the graph on the same vertex set in which two distinct vertices are adjacent if their distance in $H$ is at most $2$. In this case $H$ is called the \emph{square root} of $G$. Also, recall that the \emph{girth} of a graph is the length of its shortest cycle (or $\infty$ for a tree). The neighbourhood $N_H(u)$ of $u$ will be the set consisting of $u$ and its adjacent vertices in $H$. By $\dist_H(u,v)$ we denote the distance between two vertices in $H$. 

We investigate the uniqueness of square roots of graphs. Ross and Harary \cite{Ross} proved the following theorem:
\begin{itemize}
\item[(1)] If $T_1$ and $T_2$ are two trees such that $T_1^2$ and $T_2^2$ are isomorphic, then $T_1$ and $T_2$ are isomorphic.
\end{itemize}
This was recently improved by Farzad et al. \cite{4auth} who proved:
\begin{itemize}
\item[(2)] If $H_1$ and $H_2$ are two graphs of girth at least $7$ such that $H_1^2$ and $H_2^2$ are isomorphic, then $H_1$ and $H_2$ are isomorphic.
\end{itemize}
In the next section we prove the best possible result, which is:
\begin{itemize}
\item[(3)] If $H_1$ and $H_2$ are two graphs of girth at least $6$ such that $H_1^2$ and $H_2^2$ are isomorphic, then $H_1$ and $H_2$ are isomorphic.
\end{itemize}
The key idea behind (1) and (2) is that each maximal clique of the square corresponds to the neighbourhood of some vertex in the root. This fails in the case of roots of girth $6$. For example, the vertices $1,3,5$ of the cycle $C_6$ form a maximal clique in $C_6^2$ even though they do not induce a star in $C_6$. This is where we will need a new idea to prove (3).

\section{Proof of the theorem}
\label{section:square}

Let $H$ be a graph of girth at least $6$ on the vertex set $V$ and let $G=H^2$. We have the following easy observations:
\begin{itemize}
\item[(*)] If there is a path from $u$ to $v$ in $H$ of length exactly $3$, then $u\not\in N_G(v)$ (otherwise there would be a cycle in $H$ of length at most $5$).
\item[(**)] If $uv\in E(H)$ then $N_G(u)\cap N_G(v)=N_H(u)\cup N_H(v)$. Indeed, the inclusion $\supseteq$ is obvious. To prove $\subseteq$ note that if some vertex $w\in N_G(u)\cap N_G(v)$ was adjacent to neither $u$ nor $v$ in $H$, then it would be in distance $2$ from both of them, which would yield a 5-cycle in $H$.
\end{itemize}

We start with a lemma which can also be deduced from \cite{4auth}. The notation $H_1=H_2$ means that two graphs are equal (the same vertex set and the same edges), not just isomorphic.

\begin{lemma}
Let $H_1$ and $H_2$ be graphs of girth at least $6$ on the vertex set $V$. Suppose that $G=H_1^2=H_2^2$ and that $u,v,w\in V$ are three vertices such that $uvw$ is a path in both $H_1$ and $H_2$. Then $H_1=H_2$.
\label{lemma:easy}
\end{lemma}
\begin{proof}
Let $H$ be any graph of girth at least 6 such that $G=H^2$. The following statements follow easily from (*) and (**):

\begin{itemize} 
\item If $xyz$ is a path in $H$, then (see Fig. \ref{fig:x})
	$$N_H(x)=(N_G(x)\cap N_G(y))\setminus N_G(z)\cup\{x,y\}.$$

\begin{figure}[h!]
\includegraphics[scale=1]{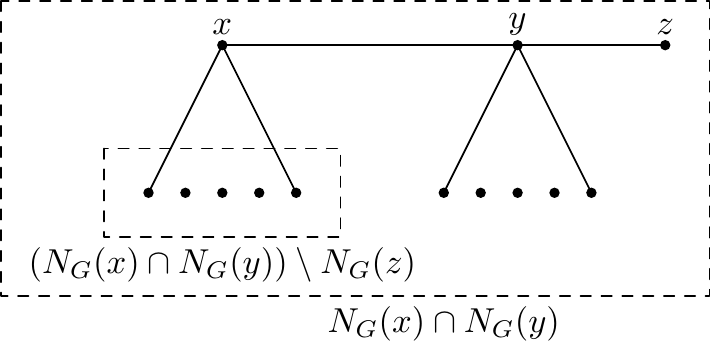}
\caption{The structure of the neighbourhoods of $x$ and $y$ in $H$.}
\label{fig:x}
\end{figure}

\item If $y$ is of degree $1$ in $H$ and $xy\in E(H)$ then  
	$$N_H(x)=N_G(y).$$
\end{itemize}
With the above formulas, given one path $uvw$ of $H$ one can recursively compute all the edges of $H$ using only the information from $G$, so the square root of $G$ with this distinguished path is unique. This ends the proof.
\end{proof}
Clearly, it suffices to prove our main result with the assumption ``$H_1^2$ and $H_2^2$ are isomorphic'' replaced by ``$H_1^2=H_2^2$''. This is what we now prove:
\begin{theorem}
\label{theorem:6iso}
Suppose $H_1$ and $H_2$ are two graphs of girth at least $6$ such that $G=H_1^2=H_2^2$. Then $H_1$ and $H_2$ are isomorphic.
\end{theorem}
\begin{proof}
Let $V$ be the common vertex set of $H_1$, $H_2$ and $G$. If $uvw$ is a path in both $H_1$ and $H_2$ for some $u,v,w$ then $H_1=H_2$ by the previous lemma. Therefore we may assume that for every $v$ the set $X_v=\{u: uv\in E(H_1)\cap E(H_2)\}$ has at most 1 element. Define the following map $f:V\longrightarrow V$:
\begin{itemize}
\item if $|X_v|=0$ then $f(v)=v$,
\item if $|X_v|=1$ then $f(v)$ is the unique element of $X_v$.
\end{itemize}
Clearly $f$ is an involution. 

We shall first prove two statements:
\begin{itemize}
\item (A) if $uv\in E(H_1)$, $|X_v|=1$ and $u\neq f(v)$ then $|X_u|=0$,
\item (B) if $uv\in E(H_1)$ and $|X_v|=0$ then $|X_u|=1$.
\end{itemize}

\textit{Proof of (A).} Let $v$ be a vertex with $|X_v|=1$ and let $f(v)=w$, meaning that $vw$ is an edge in both $H_1$ and $H_2$. Let $u$ be any neighbour of $v$ in $H_1$, other than $w$. We will show that $|X_u|=0$. Suppose, on the contrary, that $z\in X_u$ (Fig.\ref{fig:6}.). Then $\dist_{H_1}(w,z)=3$, so, by (*), $z\not\in N_G(w)$. Since $u$ and $v$ are not neighbours in $H_2$, but $u\in N_G(v)\cap N_G(w)$, the property (**) implies that $uw$ is an edge in $H_2$. However $uz$ is also an edge in $H_2$, so $z\in N_G(w)$. This contradiction proves that $|X_u|=0$ for all neighbours $u$ of $v$ in $H_1$ other than $w$.

\begin{figure}[h!]
\includegraphics[scale=0.8]{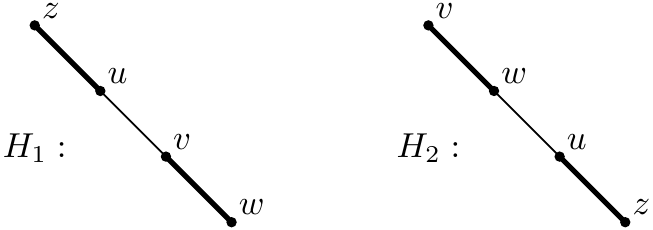}
\caption{Illustration for the proof of (A) in Theorem \ref{theorem:6iso}. The bold edges are present in both $H_1$ and $H_2$.}
\label{fig:6}
\end{figure}

\textit{Proof of (B).} Let $v$ be a vertex with $|X_v|=0$. Let $u$ be adjacent to $v$ in $H_1$. We will show that $|X_u|=1$. Suppose, on the contrary, that $|X_u|=0$. In $H_2$ the vertex $v$ must be in distance $2$ from $u$, so there is an $x$ such that $uxv$ is a path in $H_2$. In particular, $x\in N_G(v)\cap N_G(u)$, so by (**) we have $x\in N_{H_1}(u)\cup N_{H_1}(v)$. This is a contradiction since $x$ would be adjacent to either $u$ or $v$ in both $H_1$ and $H_2$, which is impossible by $X_v=X_u=\emptyset$.

\textit{Proof of the theorem.} Now we prove that $f$ (treated as a map of graphs $H_1\longrightarrow H_2$) maps edges to edges. Let $uv\in E(H_1)$. 

If $|X_u|=|X_v|=1$ then, by (A), $f(u)=v$, $f(v)=u$ and $uv$ is an edge in both graphs, so $f$ takes $uv$ to an edge $vu$ in $H_2$. 

If $|X_u|=0$ and $|X_v|=1$, then let $w=f(v)$. Since $uv\not\in E(H_2)$ and $u\in N_G(v)\cap N_G(w)$, we have from (**) that $uw\in E(H_2)$ and $f$ takes $uv\in E(H_1)$ to $f(u)f(v)=uw\in E(H_2)$.

The case $|X_u|=|X_v|=0$ is not possible by (B).

To prove that $f^{-1}$ maps edges to edges one simply inverts the roles of $H_1$ and $H_2$ in the above argument (the definition of $f$ was symmetric with respect to $H_1$ and $H_2$). Therefore $f$ is an isomorphism.
\end{proof}

\section{Remarks and modifications}
\label{section:remarks}
This result is optimal in the sense that it cannot hold for girth at least $5$ because ${K_{1,4}}^2={C_5}^2=K_5$. 

The $r$-th power $H^r$ of a graph is defined analogously, that is edges in $H^r$ correspond to pairs of vertices in distance at most $r$ in $H$. Observe that regardless of the girth restriction, there can be no analogous general result for higher graph powers, because there exist non-isomorphic trees whose $r$-th power is a complete graph for all $r\geq 3$. This and the work of \cite{Lev,Lev+} suggest that one may benefit from forbidding vertices of degree one in the root. Consider the following problem: what are the minimal values of $g_1(r)$ and $g_2(r)$ for which the following statements hold:

\begin{itemize}
\item[(1)] For any two graphs $H_1$ and $H_2$ of girth at least $g_1(r)$ with no vertices of degree one, if $H_1^r=H_2^r$ then $H_1$ and $H_2$ are isomorphic.
\item[(2)] For any two graphs $H_1$ and $H_2$ of girth at least $g_2(r)$ with no vertices of degree one, if $H_1^r=H_2^r$ then $H_1=H_2$.
\end{itemize}

\begin{figure}
\begin{tabular}{cc}
\includegraphics[scale=0.7]{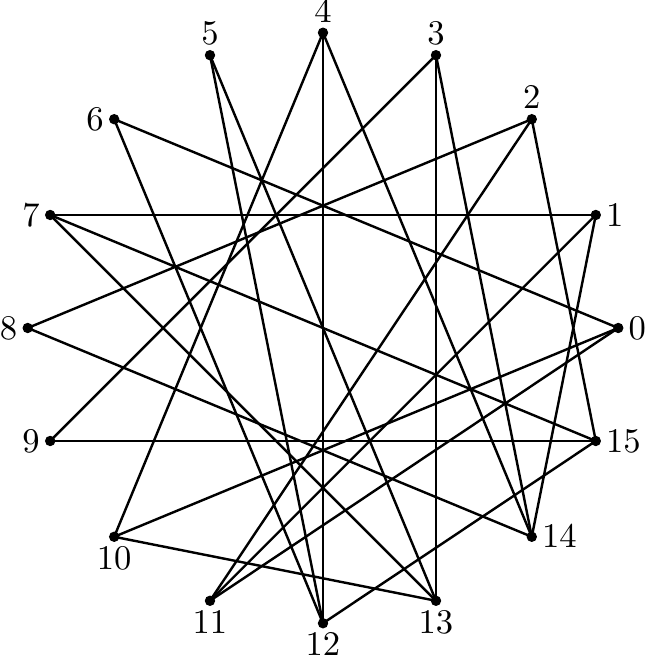} & \includegraphics[scale=0.7]{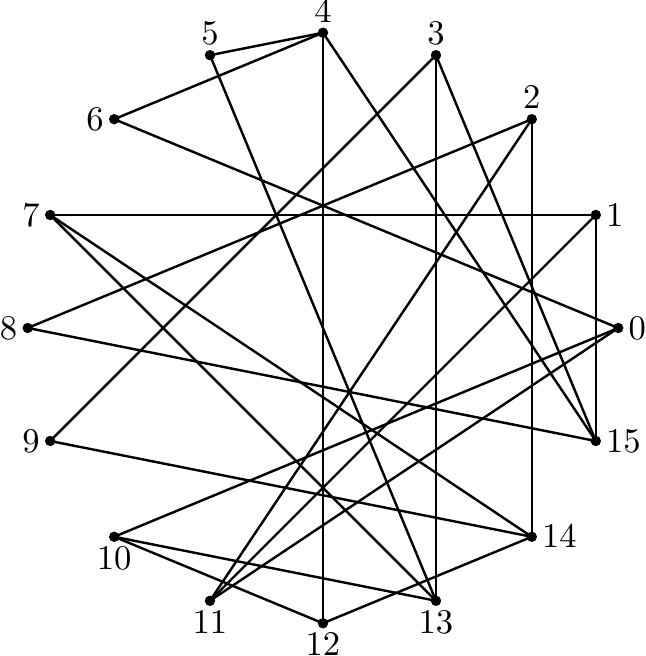}
\end{tabular}
\caption{Two non-isomorphic graphs of girth $5$, minimal vertex degree $2$ and the same square.}
\label{fig:2}
\end{figure}

For example $g_2(2)=7$, as proved in \cite{Lev+}. Our work proves that $g_1(2)\leq 6$ and this is, in fact, optimal: there exist two non-isomorphic graphs of girth $5$ and no degree one vertices having the same squares. The smallest such example is a pair of graphs on $16$ vertices shown in Fig.\ref{fig:2} (found with \cite{nauty}). Therefore $g_1(2)=6$.

It is known that $2r+3\leq g_2(r)\leq 2r+2\lceil(r-1)/4\rceil+1$ (see \cite{Lev}) and conjectured that $g_2(r)=2r+3$ for all $r$. Any nontrivial result about $g_1(r)$ (possibly in relation to $g_2(r)$) would be very interesting.

\end{document}